\documentclass[a4paper]{article}
\usepackage[latin1]{inputenc}
\usepackage[francais,english]{babel}
\usepackage{delarray,verbatim,enumerate}
\usepackage{amsmath,amsthm,amstext,amsbsy,amssymb,amsfonts,amscd,bm}

\usepackage{ifpdf}
\ifpdf
	\usepackage{aeguill}
\else
	\usepackage[T1]{fontenc}
\fi

\def\QQ{\mathbb Q}	
				
\def\Zl{\mathbb Z_\ell}		
\def\F2{\mathbb F_2}		\def\Z2{\mathbb{Z}_2}	
\def\Z2{\mathbb{Z}_2}

\def\G{\mathcal  G}		\def\F{\mathcal  F}	\def\T{\mathcal T}

\def\Gal{\operatorname{Gal}}

\def\ie{{\em i.e. }} 	\def\cf{{\em cf. }}

\newcommand{\cdast}{\mbox{\large $\circledast$}}
\def\oast{\operatorname{\circledast}}
\def\Oast{\operatornamewithlimits{\cdast}}
\newcommand{\cdplus}{\mbox{\large $\oplus$}}
\def\Oplus{\operatornamewithlimits{\cdplus}}

\newtheorem{Th}{Theorem}

\newtheorem{Cor}[Th]{Corollary}

\newtheorem{Def} [Th]{Definition}

\newtheorem{ThD}[Th]{Theorem and definition}

\begin{document}

\title{\huge\bf Note on $2$-rational fields}

\author{ Georges {\sc Gras}\; \& \; Jean-François {\sc Jaulent} }
\date{}
\maketitle
\bigskip

{\small
\noindent{\bf Résumé.} Nous déterminons le groupe de Galois de la pro-2-extension 2-ramifiée maximale d'un corps de nombres 2-rationnel.
}

\

{\small
\noindent{\bf Abstract.}  We compute the Galois group of the maximal
$2$-ramified and complexified pro-$2$-extension of any $2$-rational number field.}

\

{\small
\noindent{\bf Nota.} This short Note is motivated by the paper
``Galois $2$-extensions unramified outside $2$'' of J. Jossey
and, at this occasion, we bring into focus some classical technics of abelian
$\ell$-ramification which, unfortunately, are often ignored, especially those developped 
by  J-F. Jaulent with the $\ell$-adic class field theory, and by G. Gras in his
book on class field theory, and which considerably simplify proofs in such 
subjects; for instance, the main Theorem 2, due to J-F. Jaulent, generalizes the 
purpose of Jossey's paper in such a way.
}

\bigskip

%%%%%%%%%%%%%%%%%%%%%%%%%%%%%%%%%%%%%%%%%%%%%
%%%%%%%%%%%%%%%%%%%%%%%%%%%%%%%%%%%%%%%%%%%%%

\section{Introduction and history}

  The notions of  $\ell$-rational field and $\ell$-regular field (for a prime 
number $\ell$ and a number field $K$), independently introduced by A. Movahhedi
and T. Nguyen Quang Do in [MN], and by G. Gras and J-F. Jaulent in [GJ], 
coincide as soon as $K$ contains the maximal real subfield of the field of
$\ell$th roots of unity, thus especially for~$\ell=2$.\smallskip

\begin{itemize}

\item  The $\ell$-regularity expresses the triviality of the 
regular $\ell$-kernel of $K$ (\ie  the kernel, in
the $\ell$-part of the  universal group ${\rm K}_2(K)$, of
Hilbert symbols attached to the non-complex places not dividing $\ell$).

\item  The $\ell$-rationality 
traduces the pro-$\ell$-freeness of the Galois group
$\G_K := \Gal(M_K/K)$ of the maximal pro-$\ell$-extension 
$\ell$-ramified $\infty$-split $M_K$ of $K$
(\ie unramified at the finite places\footnote{According to the conventions 
of the $\ell$-adic class field theory (\cf 
[G$_3$, Ja]), we never speak of ramification at infinity but of
 {\it complexification} of real places.} not dividing
$\ell$ and totally split at the infinite places).

\end{itemize}
More precisely, let $c_K$ be the number of complex places of $K$; let $\bm\mu_K$ (resp. 
$\bm\mu_{K_\mathfrak l}$)\linebreak  be the $\ell$-group of roots of unity in $K$ (resp. in the localization
$K_\mathfrak l$); and let\smallskip

\centerline{$V_K := \{ x \in K^{\times} | \,  \ x \in 
K^{\times \ell}_{\mathfrak l} \  \forall  \mathfrak l \mid \ell \quad \& \quad  v_{\mathfrak p} (x) 
\equiv 0  \mod \ell  \  \ \forall {\mathfrak p} \nmid \ell\infty \}$}\smallskip

\noindent be the group of $\ell$-hyperprimary elements in $K^\times$.
Then, with these notations, from [JN, Th.1.2] or [G$_3$, IV.3.5, III.4.2.3], the
$\ell$-rationality of $K$ may be expressed as follows:

\setcounter{Th}{-1}
\begin{ThD}
The following conditions are equivalent:
 \begin{itemize}
\setlength{\itemsep}{-0,3mm}
\item[] (i) The Galois group $\mathcal G_K$ is a free pro-$\ell$-group on $1 + c_K$ 
generators.
\item[] (ii)  The abelianization $\G^{ab}_K$ of  $\mathcal G_K$ is
a free $\Zl$-module of dimension $1+c_K$.
\item[] (iii)   The field $K$ satisfies the Leopoldt conjecture (for
the prime $\ell$) and the torsion submodule $\mathcal T_K$ of
$\G^{ab}_K$ is trivial.
\item[] (iv) One has the equalities:\quad $V_K =K^{\times \ell}$ \quad \& \quad 
${\rm rk}_\ell\,( {\bm\mu}_K) = \sum_{{\mathfrak l} | \ell} {\rm 
rk}_\ell\, ({\bm\mu}_{K_{\mathfrak l}}).$
\end{itemize}
{\rm When any of these conditions is realized, the number field $K$ is said to be} $\ell$-rational.
\end{ThD}

\noindent{\bf Remark.} In case $\ell=2$, it follows from the condition {\it
(iv) } above that a 2-rational field has  {\it a single} place above 2.
\medskip

The premises of the notion of $\ell$-regularity go back to the
works of G. Gras, mainly to his note on the ${\rm K}_2$ of number fields
[G$_2$, II, \S\,2; III, \S\S\, 1, 2], 
whereas the notion of $\ell$-rationality
appears (in a hidden form) in the work of H. Miki 
[Mi] concerning the study of a sufficient condition for the 
Leopoldt conjecture, as well as those of K. Wingberg 
[W$_1$, W$_2$], concerning the same condition.\smallskip

Movahhedi's thesis and
the above papers [GJ, MN] characterised the going up
for  $\ell$-rationality
in any $\ell$-extension in terms of $\ell$-
primitivity of the ramification (a definition given in
[G$_2$, III, \S\,1] from the use of the Log function defined in
[G$_1$]), a property which was unknown in the preceding approaches.\smallskip

For instance, this gives immediately that if $K$ is a $\ell$-extension of $\QQ$,
a N.S.C. for $K$ to be $\ell$-rational
is that 
$K/\QQ$ be $\ell$-ramified, or that
$K/\QQ$ be $\{p, \ell \}$-ramified, where $p \ne \ell$
is a prime $\equiv 1 \bmod (\ell)$ such that 
$p \not\equiv \pm 1 \bmod (8)$ if
$\ell = 2$ and $p \not\equiv 1 \bmod (\ell^2)$ if $\ell\ne 2$
(\cf  [G$_3$, IV.3.5.1] giving Jossey's examples [Jo]).\smallskip

A synthesis of these results is given in [JN] and a systematic
exposition is developped in the book of G. Gras ([G$_3$, III, 
\S\,4,\,(b); IV, \S\,3,\,(b); App., \S\,2]). \smallskip

Various generalizations of these notions have been studied by  O. 
Sauzet and J-F. Jaulent (\cf [JS$_1$, JS$_2$]), especially in the 
case $\ell=2$ which is, as usual, the most triky; in particular, they 
introduce the notion of $2$-birational fields.
\medskip

Very recently, J. Jossey [Jo] introduced a notion 
of $\ell$-rationality which is incompatible with the classical one
(for  $\ell=2$, as soon as $K$ contains 
real embeddings) and is unlucky since it does not 
apply to the field of rationals  $\QQ$.\smallskip

For these reasons, to avoid any confusion, we propose to speak, in his 
context, of  $2$-{\it superrational} fields. More precisely:

\begin{Def}
 {\rm Let $K$ be a number field with  $r_K$ real places and 
 $c_K$ complex places, $M'_K$ the
 maximal $2$-ramified pro-$2$-extension  of $K$, 
and $M_K$ the maximal subextension of  $M'_K$ totally split at the
infinite places. We say that $K$ is:}
 \begin{itemize}
\setlength{\itemsep}{-0,3mm}
\item[] (i)   $2$-superrational, {\rm if  $\G'_K:=\Gal(M'_K/K)$ is pro-$2$-free;}
 
\item[] (ii)  $2$-rational, {\rm if its quotient $\G_K:=\Gal(M_K/K)$ is pro-$2$-free.}
\end{itemize}
\end{Def}

The purpose of the next section is to determine the structure of the Galois group 
$\G'_K$ when the number field $K$ is $2$-rational.\medskip

Our  proof relies on the functorial properties of $\ell$-ramification theory.

%%%%%%%%%%%%%%%%%%%%%%%%%%%%%%%%%%%%%%
%%%%%%%%%%%%%%%%%%%%%%%%%%%%%%%%%%%%%%

\section{\!\!\! Main Theorem: description of  $\G'_K = \Gal(M'_K/K)$}
 
Our very simple result has the following statement:

\begin{Th} 
Let $K$ be a $2$-rational number field having  $r_K$ 
real places and  $c_K$ complex places. The Galois group
$\G'_K:=\Gal(M'_K/K)$ of the maximal $2$-ramified
 pro-$2$-extension $M'_K$ of $K$ is the pro-$2$-free product
$$\G'_K \; \simeq \; \Z2^{\oast (1+c_K)}
\;\Oast \; (\mathbb Z /2\mathbb Z )^{\oast r_K}$$
of $\,(1+c_K)$ copies of the procyclic group $\Z2$
and of $\,r_K$ copies of $\mathbb Z/2\mathbb Z$.\smallskip
\end{Th}

\begin{Cor} The $2$-rational number fields which are  $2$-superrational
are the totally imaginary ones.
\end{Cor}

\medskip

\begin{proof} Consider the quadratic
extension  $L=K[i]$ generated by the 4th roots of unity.
It is $2$-ramified over  $K$, thus thanks to the going up theorem 
of  [GJ, MN] ({\it  cf. e.g.} [JN, Th. 3.5] or [G$_3$, IV.3.4.3, (iii)]), it is
 $2$-rational, then  $2$-superrational since it is totally imaginary.
 In other words, the Galois group $\G_L=\G'_L$ of the maximal
 $2$-ramified pro-$2$-extension $M_L$ of $L$ is pro-$2$-free.

Since the quadratic extension  $L/K$ is $2$-ramified, $M_L$ is also the 
maximal $2$-ramified pro-$2$-extension $M'_K$ of $K$; the Galois group
$\G'_K$ is {\it  potentially free} since it contains the pro-$2$-free open subgroup
$\G_L$ of index $2$ in $\G'_K$.

  As in [Jo], the results of  W. Herfort and P. 
Zalesskii (\cf [HZ, Th.0.2]) give the existence of a 
finite familly $(\F_i)_{i=0,\ldots ,k}$ of free pro-$2$-groups
on respectively $d_0,\ldots, d_k$ generators (where $k$ is the number
of conjugacy classes of subgroups of order $2$ in $\G'_K$),
such that:
$$\G'_K \; \simeq \; \F_0 \Oast \Big( \Oast_{i=1}^k (\F_i \times \mathbb Z/2\mathbb Z) 
\Big) .$$
In particular, the abelianisation $\G^{'ab}_K$ of $\G'_K$ admits
the direct decomposition:
\medskip

\centerline{$\G^{'ab}_K \; \simeq \; \Z2^{d_0} \Oplus \left( 
\Oplus_{i=1}^k (\Z2^{d_i} \oplus \mathbb Z/2\mathbb Z) \right) \; \simeq 
\; \Z2^{ d_0+d_1+ \cdots +d_k} \oplus (\mathbb Z/2\mathbb Z)^k$.}\smallskip

\noindent Since the $2$-rational field $K$ satisfies the Leopoldt conjecture, we 
get $\sum_{i=0}^k d_i = 1+c_K$ 
as well as the isomorphism $\T'_K := {\rm tor}_{\Z2}(\G^{'ab}_K) \simeq (\mathbb Z/2\mathbb Z)^k$.
Moreover
$\T_K := {\rm tor}_{\Z2}(\G^{ab}_K)= 1$, so that $\T'_K$ is generated by 
the decomposition
groups of the real places of $K$ which are deployed, a key argument 
of class field theory (\cf [Ja] or [G$_3$, 
III.4.1.5]) giving $k=r_K$.

Now the pro-$2$-decomposition of $\G'_K$ clearly shows that the minimal
number of generators $d(\G'_K)$ and of relations $r(\G'_K)$, defining
 $\G'_K$ as a pro-$2$-group, are:\medskip
 
\centerline{$d(\G'_K) = k + \sum_{i=0}^k d_i = r_K+1+c_K \hbox{\ 
and\ } r(\G'_K) = \sum_{i=1}^k (1+d_i) =  d(\G'_K)  - d_0$.}\medskip

It is well-known  ({\it cf. e.g.} [G$_3$, App., Th.2.2, (i)]) that one has\footnote{This 
argument is equivalent to the use of the formulas of 
\v Safarevi\v c ({\it  cf.} [Sa] or [NSW, Th.8.7.3])}:
$$d(\G'_K) - r(\G'_K) =\dim_{\,\mathbb F_2}(H^1(\G'_K,\mathbb F_2))
- \dim_{\,\mathbb F_2}(H^2(\G'_K,\mathbb F_2)) = 1+c_K.$$

Thus we obtain $d_0=1+c_K$, giving $d_i=0$ for $1\leq i\leq k$, 
then the expected result.
\end{proof}

%%%%%%%%%%%%%%%%%%%%%%%%%%%%%%%%%%%%%%%%
%REFERENCES
%%%%%%%%%%%%%%%%%%%%%%%%%%%%%%%%%%%%%%%%

\bigskip

\begin{tabular}{l p{1.5 cm}  l}
Georges {\sc Gras}					&{}	    & Jean-François {\sc Jaulent}\\
													&{} 	& Institut de Mathématiques\\
Villa la Gardette						&{} 	& Université de {\sc Bordeaux}\\
 chemin Ch\^ateau Gagni\`ere	&{} 	&  351, cours de la Libération\\
 F-38520 Le Bourg d'Oisans	&{}		& F-33405 TALENCE  Cedex\\
g.mn.gras@wanadoo.fr			 	&{} 	& jaulent@math.u-bordeaux1.fr
\end{tabular}

\end{document}